\documentclass[reqno,11pt]{amsart}

\usepackage{a4wide,graphicx}
\usepackage{subfigure}
\usepackage{color}
\usepackage{txfonts}
\usepackage{mathrsfs}
\usepackage[all]{xy}
\input xy
\xyoption{all}
\newtheorem{theorem}{Theorem}[section]
\newtheorem{corollary}[theorem]{Corollary}
\newtheorem{propos}[theorem]{Proposition}
\newtheorem{example}[theorem]{Example}
\newtheorem{lemma}[theorem]{Lemma}

\renewcommand{\mod}{\operatorname{mod}\nolimits}

\newcommand{\add}{\operatorname{add}\nolimits}
\newcommand{\Sub}{\operatorname{Sub}\nolimits}

\newcommand{\soc}{\operatorname{soc}\nolimits}
\newcommand{\Hom}{\operatorname{Hom}\nolimits}

\newcommand{\End}{\operatorname{End}\nolimits}

\newcommand{\Ext}{\operatorname{Ext}\nolimits}

\newcommand{\id}{{\operatorname{id}}}

\begin{document}
\title{$2$-Auslander algebras associated with reduced words in Coxeter groups}
\author{Osamu Iyama}
\address{Graduate School of Mathematics, Nagoya University}
\email{iyama@math.nagoya-u.ac.jp}
\thanks{Both authors were supported by the Storforsk-grant 167130 from the Norwegian Research Council}
\thanks{The first author was supported by JSPS Grant-in-Aid for Scientific Research 21740010 and 21340003}

\author{Idun Reiten}
\address{Insitutt for matematiske fag,
Norges Teknisk-Naturvitenskapelige Universitet,
N-7491 Trondheim, Norway}
\email{idunr@math.ntnu.no}
\maketitle

\begin{abstract}
In this paper we investigate the endomorphism algebras of standard cluster tilting objects in the stably $2$-Calabi-Yau categories $\Sub{\Lambda_w}$ with elements $w$ in Coxeter groups in \cite{BIRSc}. They are examples of the $2$-Auslander algebras introduced in \cite{I1}.
Generalizing work in \cite{GLS1} we show that they are quasihereditary, even strongly quasihereditary in the sense of \cite{R}. We also describe the cluster tilting object giving rise to the Ringel dual,
and prove that there is a duality between $\Sub{\Lambda_w}$ and the category $\mathcal{F}(\Delta)$ of good modules over the quasihereditary algebra. When $w = uv$ is a reduced word, we show that the $2$-Calabi-Yau triangulated category $\underline{\Sub}\Lambda_v$ is equivalent to a specific subfactor category of $\underline{\Sub}\Lambda_w.$
This is applied to show that a standard cluster tilting object $M$ in $\Sub{\Lambda_w}$ and the cluster tilting object $\Lambda_w\oplus\Omega{M}$ lie in the same component in the cluster tilting graph.
\end{abstract}

\section*{Introduction}
Let $A$ be a finite dimensional algebra over an algebraically closed field $k,$ where $\mbox{id}_AA\leq 1$ and $\mbox{id}$ denotes injective 
dimension.
Then the category $\Sub A$ of submodules of free $A$-modules of finite rank is an extension closed subcategory of the category $\mod A$ of finitely generated $A$-modules.
Further, $\mathscr{C} = \Sub A$ is a Frobenius category, that is, the projective and injective objects coincide and there are enough projective and enough injective objects.

In this paper we consider the important cases when $\mathscr{C}$ is stably $2$-Calabi-Yau, that is, the stable category $\underline{\mathscr{C}}$ is a $2$-Calabi-Yau triangulated category.
Let $M$ be a cluster tilting object in $\mathscr{C},$ that is, $\Ext_\mathscr{C}^1(M,M) = 0$ and if $\Ext_\mathscr{C}^1(M,X) = 0,$ then $X$ is a summand of a finite direct sum of copies of $M.$
The endomorphism algebras $\End_\mathscr{C}(M)$ belong to the class of algebras called $2$-Auslander algebras (see \cite{I2}\cite{I3}), and they are known to have global dimension at most $3$ \cite{I2}.

In this paper we deal with the finite dimensional factor algebras $\Lambda_w$ of preprojective algebras $\Lambda$ associated with elements $w$ in Coxeter groups \cite{IR}\cite{BIRSc}.
Then $\mbox{id}_{\Lambda_w}\Lambda_w\leq 1,$ and we know that $\mathscr{C}_w  = \Sub\Lambda_w$ is stably $2$-Calabi-Yau, and $\underline{\mathscr{C}_w} = \underline{\Sub}\Lambda_w$ is triangulated $2$-Calabi-Yau (see also \cite{GLS1} for the case of adaptable words).
We consider mainly cluster tilting objects $M$ which are associated with reduced expressions of $w,$ which are called standard cluster tilting objects. 
In this case we show that the $2$-Auslander algebras $\End_{\Lambda_w}(M)$ are quasihereditary, even strongly quasihereditary in the terminology of \cite{R}, and that the Ringel dual quasihereditary algebra of $\End_{\Lambda_w}(M)$ is $\End_{\Lambda_w}(\widetilde{\Omega}M),$
where $\widetilde{\Omega}M$ denotes the direct sum of $\Lambda_w$ and the syzygy module $\Omega M$ of $M$.

Using mutation of cluster tilting objects in $\Sub\Lambda_w$ or in $\underline{\Sub}\Lambda_w,$ there is an associated graph called the cluster tilting graph, where the vertices correspond to the isomorphism classes of basic cluster tilting objects.
It is an important open problem whether this graph is connected. This is known to be the case for cluster categories of finite dimensional hereditary algebras (see \cite{BMRRT}), and for cluster categories of coherent sheaves on weighted projective lines in the tubular case \cite{BKL}.
Here we show that if $M$ is a standard cluster tilting object in $\Sub\Lambda_w,$ then $\widetilde{\Omega}M,$ which gives rise to the Ringel dual of $\End_{\Lambda_w}(M),$ lies in the same component as $M.$
In order to prove this, we construct, for a reduced word $w = uv,$ an embedding of $\mathscr{C}_v$ into $\mathscr{C}_w,$ which induces an equivalence of triangulated categories from $\underline{\mathscr{C}_v}$
to a subfactor triangulated category of $\underline{\mathscr{C}_w},$ using a general construction investigated in \cite{IY}. This is also of interest in its own right.

The paper is organised as follows. In Section \ref{sc1} we give some background material on Hom-finite $2$-Calabi-Yau categories with cluster tilting objects, in particular we deal with those associated with elements in Coxeter groups.
We also give basic definitions and facts about $2$-Auslander algebras and about (strongly) quasihereditary algebras.
In Section \ref{sc2} we show that the endomorphism algebras of standard cluster tilting objects $M$ in $\Sub\Lambda_w$ are strongly quasihereditary, with Ringel dual given  by $\widetilde{\Omega}M.$
We also show  a strong relationship between the category $\Sub\Lambda_w$ and the category $\mathcal{F}(\Delta)$ of $\End_{\Lambda_w}(M)$-modules with good  filtrations (see \cite{R2}).  
In Section \ref{sc3} we discuss the embedding of $\mathscr{C}_v = \Sub\Lambda_v$ into $\mathscr{C}_w = \Sub\Lambda_w,$ which induces our desired equivalence of triangulated categories. 
Then we apply this in Section \ref{sc4} to show that $M$ and $\widetilde{\Omega}M$ lie in the same component in the cluster tilting graph.

This work was inspired by the work on quasihereditary algebras for the case of adaptable words in \cite{GLS1}, and was presented in Mexico City (December 2008), Bielefeld (June 2009) and Durham (July 2009). 
Most of the results in Sections \ref{sc2} and \ref{sc4} have later also been  proved using different methods in \cite{GLS2}. A further generalization of our class of quasihereditary algebras in Section \ref{sc2} has been announced in \cite{S}.

\medskip\noindent{\bf Acknowledgements }
Part of this work was done while the authors participated in the conference ``The Representation Dimension of Artin Algebras'' (Bielefeld, May 2008)
and the first author visited NTNU during March and August 2009. He would like to thank the people in Bielefeld and Trondheim for their hospitality.

\section{Background}\label{sc1}
Throughout this paper all modules are left modules, and the composition $fg$ of morphisms means first $g$, then $f$.
In this section we give some background material on $2$-Calabi-Yau ($2$-CY for short) categories, $2$-Auslander algebras and quasihereditary algebras.

\subsection{$2$-CY categories}
Let $A$ be a finite dimensional basic $k$-algebra. An extension closed subcategory $\mathscr{C}$ of $\mod A$ is called \emph{Frobenius} if the projective and injective objects coincide, and there are enough projective and enough injective objects. 
Then $\mathscr{C}$ is \textit{stably $2$-CY} if the stable category $\underline{\mathscr{C}}$ is a $2$-CY triangulated category. 
An object $M$ in $\mathscr{C}$ (or $\underline{\mathscr{C}}$) is \textit{cluster tilting} if $\Ext_\mathscr{C}^1(M,X) = 0$ if and only if $X$ is in $\mbox{add}M.$ 
Here $\mbox{add}M$ denotes the full additive subcategory of $\mathscr{C}$ (or $\underline{\mathscr{C}}$) whose objects are finite direct sums of copies of $M.$

When $\mathscr{C}$ is  Hom-finite triangulated $2$-CY, there is a way of constructing subfactors of $\mathscr{C}$ which are again Hom-finite triangulated $2$-CY \cite{IY} (see also \cite{BIRSc}). 
Let  $D$ be a rigid object in $\mathscr{C}$, that is $\Ext_\mathscr{C}^1(D,D) = 0,$ and consider  $D^{\bot_1} = \{X\in\mathscr{C}; \Ext_\mathscr{C}^1( D,X) = 0\}.$ 
Then the factor  category $D^{\bot_1}/{\mbox{add}D}$ is triangulated $2$-CY, and there is a one-one correspondence between the cluster tilting objects in $\mathscr{C}$ containing  $D$ as a summand, and the cluster tilting objects in $ D^{\bot_1}/{\mbox{add}D}.$

Let $M = M_1\oplus\ldots\oplus M_n$ be a cluster tilting object in the stably $2$-CY category $\mathscr{C},$ where the $M_i$ are indecomposable and nonisomorphic. 
Assume that $M_i$ is not projective for $1\leq i\leq m$ and $M_i$ is projective for  $m < i\leq n.$  For each $i = 1,\ldots,m$ there is a unique indecomposable object $M_i^*\nsimeq M_i$ such that $\mu_i(M) = (M/M_i)\oplus M_i^*$ is a cluster tilting object in $\mathscr{C}.$ 
This gives rise to a graph, the \textit{cluster tilting graph}, where the vertices correspond to cluster tilting objects up to isomorphism. For each $M$ there are $m$ vertices $\mu_1(M),\ldots,\mu_m(M)$ connected to $M$ by an edge.

\subsection{$2$-CY categories associated with words in Coxeter groups}
An important class of (stably) $2$-CY categories are those associated with reduced  words in Coxeter groups \cite{IR}\cite{BIRSc}. 
Let $Q$ be a finite connected quiver with $n$ vertices and no oriented cycles, and let $W_Q$ be the associated Coxeter group, with generators $s_1,\ldots,s_n,$ and let $\Lambda$ be the associated preprojective algebra. 
For each $i = 1,\ldots,n,$ let $I_i := \Lambda(1-e_i)\Lambda,$ where $e_i$ is the idempotent element at the vertex $i$ of $Q.$  Let $w = s_{i_1}\ldots s_{i_t}$ be a reduced word in $W_Q.$ 
Then the ideal $I_w := I_{i_1}\ldots I_{i_t}$ is independent of the reduced expression of $w,$ and if $Q$ is non-Dynkin, $I_w$ is a tilting $\Lambda$-module of projective dimension at most one. 
The algebra $\Lambda_w := \Lambda/I_w$ is a finite dimensional algebra with $\mbox{id}_{\Lambda_w}\Lambda_w\leq 1,$ so that $\Sub\Lambda_w$ is a Frobenius category. 
Further $\Sub\Lambda_w$ is stably $2$-CY, and the stable category $\underline{\Sub}\Lambda_w$ is triangulated $2$-CY.

We write $\mathbf{w}$ when we mean the reduced expression of $w.$
For $\mathbf{w} = s_{i_1}\ldots s_{i_t},$ let
\begin{equation}\label{standard cluster tilting}
M_j=P_{i_j}/(I_{i_1}\ldots I_{i_j})P_{i_j}\ \mbox{ and }\ M_\mathbf{w} = M_1\oplus\ldots\oplus M_t,
\end{equation}
where $\Lambda = P_1\oplus\ldots\oplus P_n,$ and $P_i$ is the indecomposable projective $\Lambda$-module associated with the vertex $i.$ 
Then $M_\mathbf{w}$ is a cluster tilting object in $\Sub\Lambda_w$ and in $\underline{\Sub}\Lambda_w,$ which we call a \emph{standard cluster tilting object}. 
There is clearly only a finite number of standard cluster tilting objects in $\Sub\Lambda_w$, and they are all known to lie in the same component of the cluster tilting graph \cite{BIRSc}.
For a vertex $i$ in the quiver $Q$, let $i_{l_1},\ldots,i_{l_t}$ be the ordered vertices of type $i$ in $\mathbf{w}$.
Then we have the epimorphisms $M_{l_j}\rightarrow M_{l_{j-1}}$, and we denote the kernels by $L_{l_j}$ and call them \emph{layers} as in \cite{AIRT}.

\subsection{$2$-Auslander algebras}
The \emph{Auslander algebras} are finite dimensional algebras which by definition are the endomorphism algebras $\End_A(M)$ when $M$ is an additive generator of $\mod A$ for a finite dimensional algebra $A$ of finite representation type \cite{A}. 
They are characterized as being algebras $\Gamma$ of global dimension at most two and dominant dimension at least two, that is, in the minimal injective resolution $0\rightarrow\Gamma\rightarrow I_0\rightarrow I_1\rightarrow I_2\rightarrow 0$ of $\Gamma,$ both $I_0$ and $I_1$ are projective.

In \cite{I2}\cite{I3} the more general concept of $n$-Auslander algebras was introduced for $n\geq 1,$ where the $1$-Auslander algebras are the Auslander algebras above. 
For $n = 2$ we have the following. Let $U$ be a cotilting module and ${}^\bot U = \{X\in\mod A; \mbox{ Ext}_A^i(X, U) = 0\ (i>0)\}.$ 
Then $M$ is a cluster tilting object in $ {}^\bot U$ if $\mbox{add}M = \{X\in{}^\bot U; \mbox{ Ext}_A^1(X,M) = 0\} = \{Y\in{}^\bot U; \mbox{ Ext}_A^1(M,Y) = 0\}.$ 
When $\id{}_AU\le 1$, the algebras $\End_A(M)$ are ($1$-relative) \emph{$2$-Auslander algebras}.  
They are the finite dimensional algebras $\Gamma$ with $\mbox{gl.dim}\ \Gamma\leq 3,$ and if $0\rightarrow\Gamma\rightarrow I_0\rightarrow I_1\rightarrow I_2\rightarrow I_3\rightarrow 0$ is a minimal injective resolution of $\Gamma,$ then $I_0, I_1$ and $I_2$ have projective dimension at most $1.$ 
When $A = \Lambda_w$ and $ U = \Lambda_w$ for an element $w$ in a Coxeter group and $M$ is a standard cluster tilting object in $\Sub\Lambda_w = {}^\bot\Lambda_w,$ the algebras $\End_A(M)$ are $2$-Auslander algebras. 
When $M$ is a standard cluster tilting object in $\Sub\Lambda_w,$ there is an explicit description of the $2$-Auslander algebras $\End_{\Lambda_w}(M)$ in terms of quivers with relations \cite{BIRSc}\cite{BIRSm}.

\subsection{Strongly quasihereditary algebras}
Let $A$ be a finite dimensional $k$-algebra. We say that an ideal $I$ of $A$ is \emph{heredity} if $I^2=I$, $I$ is a projective $A$-module and $\End_A(I)$ is a semisimple algebra. 
We say that $A$ is \emph{quasihereditary} if there exists a chain $A\supset I_1\supset\cdots\supset I_n=0$ of ideals of $A$ such that $I_{i-1}/I_i$ is a heredity ideal of $A/I_i$ for any $i$ (see \cite{DR}).
This is equivalent to the following condition: Let $P_1,\ldots,P_n$ be nonisomorphic indecomposable projective $A$-modules.
For each $i = 1,\ldots,n,$ denote by $\Delta_i$ the largest factor of $P_i$ with composition factors amongst the simple modules $S_1,\ldots,S_i,$ where $S_j$ is associated to $P_j.$ 
Then $A$ is quasihereditary (with respect to the ordering $P_1,\ldots,P_n$) if and only if each $P_i$ has a filtration by the modules $\Delta_1,\ldots,\Delta_n,$ and $A$ has finite global dimension (see \cite{R2}). 
The algebra $A$ is said to be (left) \textit{strongly quasihereditary} (see \cite{R}) if each $\Delta_i$ has projective dimension at most one. 
The subcategory $\mathcal{F}(\Delta)$ of $\mod A,$ whose objects have a filtration using $\Delta_1,\ldots,\Delta_n,$ is contravariantly finite and \textit{resolving}, that is extension closed, closed under kernels of epimorphisms and contains the projectives. 
There is a cotilting module $U$ associated with $\mathcal{F}(\Delta),$ which is also a tilting module, given by the indecomposable Ext-injective modules in $\mathcal{F}(\Delta).$ 
Then we have $\mathcal{F}(\Delta) = {}^\bot{U}$, and $U$ is said to be a \emph{characteristic tilting module}. 
The algebra $\End_A(U)$ is again quasihereditary, and is called the \emph{Ringel dual} of $A$ (see \cite{R2}).
We say that $A$ is \emph{$\Delta$-serial} if the indecomposable projective $A$-modules have a unique $\Delta$-composition series. 

\section{Construction of quasihereditary algebras with additional properties}\label{sc2}
Throughout this section, let $\Lambda$ be a preprojective algebra of a finite quiver $Q$ without oriented cycles,
$w$ be an element in the Coxeter group, and $\mathbf{w}$ be a reduced expression of $w$.
We have a standard cluster tilting object $M=M_{\mathbf{w}}$ in $\Sub\Lambda_w$ and a $2$-Auslander algebra $\Gamma:=\End_{\Lambda_w}(M)$.
We show that  $\Gamma$ is strongly quasihereditary and $\Delta$-serial, and that its Ringel dual is the $2$-Auslander algebra $\End_{\Lambda_w}(\widetilde{\Omega}M)$ given by the cluster tilting object $\widetilde{\Omega}M$ in $\Sub\Lambda_w$.
We give two different approaches to proving that $\Gamma$ is quasihereditary, where the second one depends heavily on \cite{BIRSc}, while the first one is  more direct. 

With the previous notation we have the following.
\begin{theorem}
$\Gamma=\End_{\Lambda_w}(M)$ is a quasihereditary algebra.
\end{theorem}

\begin{proof}
We denote by $e$ the idempotent of $\Gamma$ corresponding to the simple direct summand $M_1=P_{i_1}/I_{i_1}P_{i_1}$ of $M$.
It suffices to show the following assertions.
\begin{itemize}
\item[(i)] $\Gamma e\Gamma$ is a heredity ideal of $\Gamma$.
\item[(ii)] $\Gamma/\Gamma e\Gamma$ is isomorphic to $\Gamma':=\End_{\Lambda_{w'}}(M')$ for the cluster tilting object $M'=M_{\mathbf{w}'}$
in $\Sub\Lambda_{w'}$ associated to the reduced expression $\mathbf{w}'=s_{i_2}\cdots s_{i_t}$.
\end{itemize}
Then the theorem follows inductively.
The statement (i) follows from the following general observation.

\begin{lemma}
Let $\Lambda$ be an algebra and $M$ be a finite length $\Lambda$-module with a simple direct summand $S$.
Let $\Gamma=\End_\Lambda(M)$ and $e$ the idempotent of $\Gamma$ corresponding to the direct summand $S$ of $M$.
Then $I=\Gamma e\Gamma$ is a heredity ideal of $\Gamma$.
\end{lemma}

\begin{proof}
We denote by $\soc_S(M)$ the sum of the simple submodules of $M$ which are isomorphic to $S$.
Then the inclusion map $f:\soc_S(M)\to M$ induces an isomorphism
$\Hom_\Lambda(M,\soc_S(M))\simeq I$ of $\Gamma$-modules.
In particular $I$ is a projective $\Gamma$-module.
Dually $I$ is a projective $\Gamma^{\rm op}$-module.

Since $\End_\Gamma(\Gamma e)$ is Morita equivalent to $e\Gamma e\simeq\End_\Lambda(S),$ which is a division algebra,
we have that $I$ is a heredity ideal.
\end{proof}

Now we show (ii).
We have a functor $F:\mod\Lambda\to\mod\Lambda$ given by
\[F(X):=X/\soc_{M_1}(X).\]
Clearly we have that $F(M)$ is $M'$ above.
Thus we have an algebra homomorphism
\[\phi=F_{M,M}:\Gamma=\End_{\Lambda}(M)\to\Gamma'=\End_{\Lambda}(M').\]
We will show that $\phi$ induces an isomorphism $\Gamma/\Gamma e\Gamma\simeq\Gamma'$.
Clearly $f\in\Gamma$ satisfies $\phi(f)=0$ if and only if $f$ factors through $\add M_1$ if and only if $f\in\Gamma e\Gamma$.

We only have to show that $\phi$ is surjective.
Fix any $g\in\End_{\Lambda}(M')$ and consider the exact sequence
\[0\to\soc_{M_1}(M)\to M\xrightarrow{p} M'\to0.\]
By \eqref{standard cluster tilting} the syzygy $\Omega_\Lambda M$ of the $\Lambda$-module $M$ satisfies $\Hom_{\Lambda}(\Omega_\Lambda M,M_1)=0$,
so we have $\Ext^1_{\Lambda}(M,M_1)=0$.
Hence the map $gp:M\to M'$ factors through $p$,
and there exists $f\in\End_{\Lambda}(M)$ such that $pf=gp$.
Then $f$ satisfies $g=F(f)=\phi(f)$, and we have the assertion.
\end{proof}
 For the second proof we first give a sufficient condition for $\End_\mathscr{C}(M)$ to be strongly quasihereditary, for an object $M$ in some additive category $\mathscr{C}$.
 %%%%

\begin{lemma}\label{Thm2_1}
Let $\mathscr{C}$ be a Hom-finite extension closed subcategory of an abelian $k$-category, and $M$ an object in $\mathscr{C}$ with $\Ext_\mathscr{C}^1(M,M) = 0.$ 
Let $\Gamma = \End_\mathscr{C}(M).$ Write $M = M_1\oplus\ldots\oplus M_n,$ where the $M_i$ are indecomposable and nonisomorphic. 
Assume that for each $M_i$ the minimal left $\emph{add}(\underset{j<i}{\bigoplus}M_j)$-approximation $f_i:M_i\rightarrow M_i'$ is surjective, and that $\emph{gl.dim}\ \Gamma<\infty.$ Then we have the following.
\begin{enumerate}
	\item $\Gamma$ is (left) strongly quasihereditary with respect to the ordering $\Hom_\mathscr{C}(M_n,M),\ldots, \Hom_\mathscr{C}(M_1,M)$ of the nonisomorphic indecomposable projective $\Gamma$-modules, and the associated factor modules of the $\Hom_\mathscr{C}(M_i,M)$ are the $\Delta_i = \Hom_\mathscr{C}(\emph{Ker}f_i,M),$ for $i = 1,\ldots,n.$
	\item If $M_i'$ is indecomposable or zero for each $i = 1,\ldots,n,$ then $\Gamma$ is $\Delta$-serial.
\end{enumerate}
\end{lemma}

\begin{proof}
	\textbf{(a)}\hspace{0.1cm} For each $i = 1,\ldots,n,$ consider the exact sequence $0\rightarrow L_i\rightarrow M_i\stackrel{f_i}{\rightarrow}M_i'\rightarrow 0,$ where $L_i = \mbox{Ker}f_i.$ 
Since $\Ext_\mathscr{C}^1(M,M) = 0,$ we have an exact sequence $0\rightarrow\Hom_{\mathscr{C}}(M_i',M)\rightarrow\Hom_{\mathscr{C}}(M_i,M)\rightarrow\Hom_{\mathscr{C}}(L_i,M)\rightarrow 0,$ and hence $\mbox{pd}_\Gamma\Hom_{\mathscr{C}}(L_i,M)\leq 1.$
	
Denote by $S_i$ the simple $\Gamma$-module which is the top of the indecomposable projective $\Gamma$-module $\Hom_{\mathscr{C}}(M_i,M).$ Assume that $S_r$ is a composition factor of $\Hom_{\mathscr{C}}(L_i,M),$ for some $r = 1,\dots,n.$ 
Then we have a nonzero map $\Hom_{\mathscr{C}}(M_r,M)\stackrel{(g,M)}{\rightarrow}\Hom_{\mathscr{C}}(M_i,M)$ such that the composition $\Hom_{\mathscr{C}}(M_r,M)\rightarrow\Hom_{\mathscr{C}}(M_i,M)\rightarrow\Hom_{\mathscr{C}}(L_i,M)$ is nonzero. 
Hence the map $\Hom_{\mathscr{C}}(M_r,M)\rightarrow\Hom_{\mathscr{C}}(M_i,M)$ does not factor through the map $\Hom_{\mathscr{C}}(M_i',M)\stackrel{(f_i,M)}{\rightarrow}\Hom_{\mathscr{C}}(M_i,M).$ 
Then $g:M_i\rightarrow M_r$ does not factor through $f_i:M_i\rightarrow M_i'.$ Since $f_i:M_i\rightarrow M_i'$ is a minimal left $\mbox{add}(\underset{j<i}{\oplus}M_j)$-approximation, it follows that $M_r$ is not a summand of $\underset{j<i}{\oplus}M_j.$ 
Hence we have $r\geq i$ for all simple composition factors $S_r$ of $\Hom_{\mathscr{C}}(L_i,M).$

We want to show that $\Hom_{\mathscr{C}}(L_i,M)$ is the largest factor of $\Hom_{\mathscr{C}}(M_i,M)$ with composition factors amongst $S_r$ for $r\geq i.$ 
This follows from the exact sequence $0\rightarrow\Hom_{\mathscr{C}}(M_i',M)\rightarrow\Hom_{\mathscr{C}}(M_i,M)\rightarrow\Hom_{\mathscr{C}}(L_i,M)\rightarrow 0,$ since the top of $\Hom_{\mathscr{C}}(M_i',M)$ only has composition factors $S_j$ with $j<i.$ 
Then we see that $\Delta_i = \Hom_{\mathscr{C}}(L_i,M),$ and hence $\mbox{pd}_\Gamma{\Delta_i}\leq 1.$

We show that each indecomposable projective $\Gamma$-module $\Hom_{\mathscr{C}}(M_i,M)$ has a $\Delta$-filtration by induction on the length. If the length of $\Hom_{\mathscr{C}}(M_i,M)$ is smallest possible, then $M_i' = 0,$ so that $\Hom_{\mathscr{C}}(M_i,M)=\Delta_i.$ 
The rest follows easily. This finishes the proof of (a).
	
\textbf{(b)}\hspace{0.1cm} The top of the $\Delta$-filtration of the indecomposable projective $\Gamma$-module $\Hom_{\mathscr{C}}(M_i,M)$ has to be $\Delta_i$, and the remaining part is $\Hom_{\mathscr{C}}(M'_i,M)$ which is zero or an indecomposable projective $\Gamma$-module.
Thus the assertion follows by induction on the length of the indecomposable projectives.
\end{proof}

We have the following direct consequence of Lemma \ref{Thm2_1}.

\begin{theorem}\label{cor2_4}
Let $M$ be a standard cluster tilting object in $\Sub\Lambda_w$ associated with a reduced expression $\mathbf{w}=s_{l_1}\ldots s_{l_t}$
of an element in a Coxeter group. Then: 
\begin{enumerate}
\item $\Gamma = \End_{\Lambda_w}(M)$ is (left) strongly quasihereditary.
\item We have $\Delta_i=\Hom_{\Lambda_w}(L_i,M)$.
\item The indecomposable projective $\Gamma$-modules are $\Delta$-uniserial.
\end{enumerate}
\end{theorem}

\begin{proof}
Since $\Gamma$ is a 2-Auslander algebra, we have $\mbox{gl.dim}\ \Gamma\le 3$.

Fix a vertex $i$ in the quiver $Q.$ Let $i_{l_1},\ldots,i_{l_t}$ be the ordered vertices of type $i$ in $\mathbf{w}$. 
Then we have a sequence of irreducible epimorphisms $M_{l_t}\rightarrow\ldots\rightarrow M_{l_1}$ in $\mbox{add}M,$ which correspond to arrows going from right to left in the quiver of $\mbox{add}M$ (see \cite{BIRSc}). 
All other arrows in the quiver go from left to right. 
By Lemma \ref{Thm2_1}, we only have to show that the epimorphisms $M_{l_j}\rightarrow M_{l_{j-1}}$ for $j\geq 2$ and $M_{l_1}\rightarrow 0$ are minimal left $\mbox{add}(\underset{r<l_j}{\oplus}M_r)$-approximations of $M_{l_j}.$   

This is easy to see directly, or one can use an idea from \cite[Th.6.6]{BIRSm}:
We consider a nonzero path $C:M_{l_j}\rightarrow M_r$ with $r<l_j.$  
On the basis of the right turning points we define $\alpha(C),$ and show as in  \cite[Th.6.6]{BIRSm} that we can replace $C$ by another path representing the same element, but with smaller $\alpha$-value. 
Then we assume that we have made a choice with $\alpha(C)$ minimal. Assume the start is not $M_{l_j}\rightarrow M_{l_{j-1}}.$ Then we have to go to the right from $M_{l_{j-1}},$ and hence have a right turning point. 
So we can reduce the $\alpha$-value and get a contradiction to the minimality, and we are done. 
\end{proof}

Our next aim is to show that $\Hom_\Lambda(\widetilde{\Omega}M,M)$ is the characteristic tilting module for the quasihereditary algebra $\Gamma.$ Since the $\Delta_i$ have projective dimension at most one, and the characteristic tilting module is filtered by the $\Delta_i,$ we know that it must have projective dimension at most one. We recall the following information.

\begin{propos}\label{U is tilting}
Let $M$ in $\Sub\Lambda_w$ be a standard cluster tilting object as before, and $\Gamma = \End_{\Lambda_w}(M)$.
\begin{enumerate}
	\item Then $\widetilde{\Omega}M$ is a cluster tilting object in $\Sub\Lambda_w,$ and $\underline{\End}_{\Lambda_w}(M)\simeq\underline{\End}_{\Lambda_w}(\Omega{M}).$
	\item $U = \Hom_{\Lambda_w}(\widetilde{\Omega}M,M)$ is a tilting $\Gamma$-module of projective dimension at most one  such that $\End_\Gamma(U)\simeq\End_{\Lambda_w}(\widetilde{\Omega}M).$
\end{enumerate}
\end{propos}

\begin{proof}
	\textbf{(a)}\hspace{0.1cm} This follows since $\mbox{id}_{\Lambda_w}\Lambda_w\leq 1,$ and hence $\Omega:\underline{\Sub}\Lambda_w\rightarrow\underline{\Sub}\Lambda_w$ is an equivalence.
	
\textbf{(b)}\hspace{0.1cm} This is \cite[Th.5.3.2]{I2}.
\end{proof}

In order to show that $ U = \Hom_{\Lambda_w}(\widetilde{\Omega}M,M)$ is the characteristic tilting module for $\Gamma,$ it will be sufficient to prove that $ U$ is in $\mathcal{F}(\Delta)$ and that $\mathcal{F}(\Delta)$ is in ${}^\bot U.$ 
For the second statement it is sufficient to show that $\Hom_{\Lambda_w}(\Sub\Lambda_w,M)$ is contained in ${}^\bot U,$ since we have already seen that the $\Delta_i$ are contained in $\Hom_{\Lambda_w}(\Sub\Lambda_w,M)$ and hence $\mathcal{F}(\Delta)\subset\Hom_{\Lambda_w}(\Sub\Lambda_w,M)$

\begin{propos}\label{prep2_4}
For any $X$ in $\Sub\Lambda_w$ we have the following.
\begin{enumerate}
	\item $\emph{pd}_\Gamma\Hom_{\Lambda_w}(X,M)\leq 1.$
	\item $\Ext^i_\Gamma(\Hom_{\Lambda_w}(X,M), U) = 0$ for all $i>0,$ so $\Hom_{\Lambda_w}(X,M)$ is in ${}^\bot U.$
\end{enumerate}
\end{propos}

\begin{proof}
	\textbf{(a)}\hspace{0.1cm} Since $M$ is a cluster tilting object in $\Sub\Lambda_w,$ we have an exact sequence
\begin{equation}\label{M-resolution}
0\rightarrow X\rightarrow M_0\stackrel{p}{\rightarrow} M_1\rightarrow 0
\end{equation}
in $\Sub\Lambda_w$ with $M_0$ and $M_1$ in $\mbox{add}M$ (see \cite{I1,BMR,KR}). We apply $\Hom_{\Lambda_w}(\,\,,M)$ to get the exact sequence $0\rightarrow\Hom_{\Lambda_w}(M_1,M)\rightarrow\Hom_{\Lambda_w}(M_0,M)\rightarrow\Hom_{\Lambda_w}(X,M)\rightarrow 0,$ showing $\mbox{pd}_\Gamma\Hom_{\Lambda_w}(X,M)\leq 1.$

\textbf{(b)}\hspace{0.1cm} Apply $\Hom_\Gamma(\,\,,\Hom_{\Lambda_w}(\widetilde{\Omega}M,M))$ to get the first exact sequence in the following diagram:
\[\xymatrix@C0.3cm@R0.3cm{
\Hom_\Gamma((M_0,M),(\widetilde{\Omega}M,M))\ar[r]\ar^\wr[d]&\Hom_\Gamma((M_1,M),(\widetilde{\Omega}M,M))\ar[r]\ar^\wr[d]&
\Ext^1_\Gamma((X,M),(\widetilde{\Omega}M,M))\ar[r]&0\\
\Hom_{\Lambda_w}(\widetilde{\Omega}M,M_0)\ar[r]&\Hom_{\Lambda_w}(\widetilde{\Omega}M,M_1)
%\ar[r]&\Ext^1_{\Lambda_w}(\widetilde{\Omega}M,X)\ar[r]&\Ext^1_{\Lambda_w}(\widetilde{\Omega}M,M_0)
}\]
The second exact sequence is obtained by applying $\Hom_{\Lambda_w}(\widetilde{\Omega}M,\,\,)$ to the exact sequence \eqref{M-resolution}, and the two isomorphisms follow since $M_0$ and $M_1$ are in $\mbox{add}M.$ 

Fix any $f\in\Hom_{\Lambda_w}(\widetilde{\Omega}M,M_1)$.
Since $\underline{\Hom}_{\Lambda_w}(\widetilde{\Omega}M,M_1)=\Ext^1_{\Lambda_w}(M,M_1)=0$, we have that $f$ factors through a projective $\Lambda_w$-module $P$.
Since $p$ in \eqref{M-resolution} is surjective, we have that $f$ factors through $p$.
%Since $\underline{\Sub}\Lambda_w$ is $2$-CY we have the commutative diagram
%\[\xymatrix@C0.3cm@R0.3cm{
%\Ext^1_{\Lambda_w}(\widetilde{\Omega}M,X)\ar[r]\ar^\wr[d]&
%\Ext^1_{\Lambda_w}(\widetilde{\Omega}M,M_0)\ar^\wr[d]&\\
%\underline{\Hom}_{\Lambda_w}(M,X[2])\ar[r]\ar^\wr[d]&
%\underline{\Hom}_{\Lambda_w}(M,M_0[2])\ar^\wr[d]\\
%\Hom_k\left(\underline{\Hom}_{\Lambda_w}(X,M),\,k\right)\ar[r]&
%\Hom_k\left(\underline{\Hom}_{\Lambda_w}(M_0,M),\,k\right)
%}\]
%Since we have seen that $\Hom_{\Lambda_w}(M_0,M)\rightarrow\Hom_{\Lambda_w}(X,M)$ is surjective, it follows directly that also $\Hom_{\underline{\Sub}\Lambda_w}(M_0,M)\rightarrow\Hom_{\underline{\Sub}\Lambda_w}(X,M)$ is surjective. 
%Then $\Ext^1_{\Lambda_w}(\widetilde{\Omega}M,X)\rightarrow\Ext^1_{\Lambda_w}(\widetilde{\Omega}M,M_0)$ is injective, and
\[\xymatrix@C1cm@R.3cm{
P\ar[dr]\ar@{.>}[d]&\widetilde{\Omega}M\ar^f[d]\ar[l]\\
M_0\ar_p[r]&M_1
}\]
Consequently, the above map $\Hom_{\Lambda_w}(\widetilde{\Omega}M,M_0)\rightarrow\Hom_{\Lambda_w}(\widetilde{\Omega}M,M_1)$ is surjective. 
Thus $\Ext_\Gamma^1(\Hom_{\Lambda_w}(X,M), U) = 0,$ and we are done by \textbf{(a)}.
\end{proof}

We can now show the desired property for $ U.$

\begin{theorem}\label{characteristic}
With the previous notation we have the following:
\begin{enumerate}
\item $ U = \Hom_{\Lambda_w}(\widetilde{\Omega}M,M)$ is in $\mathcal{F}(\Delta)$.
\item $\mathcal{F}(\Delta) = {}^\bot U.$ 
\item $ U$ is the characteristic tilting module.
\end{enumerate}
\end{theorem}

\begin{proof}
\textbf{(a)}\hspace{0.1cm} For a fixed vertex $k$ in the quiver $Q,$ we consider as before the ordered vertices $i_{l_1},\ldots,i_{l_t}$ of type $k.$ 
Then we have exact sequences $0\rightarrow L_{l_j}\rightarrow M_{l_j}\rightarrow M_{l_{j-1}}\rightarrow 0$ for $j\geq 2,$ and we have $L_{l_1} = M_{l_1}.$ The $M_{l_j}$ are  factors of the indecomposable projective $\Lambda_w$-module $P_k.$ Hence we have the exact commutative diagram
\[\xymatrix@C0.3cm@R0.3cm{
&&0\ar[d]&0\ar[d]\\
&&\Omega{M_{l_j}}\ar[r]\ar[d]&\Omega{M_{l_{j-1}}}\ar[d]\\
&&P_k\ar^\sim[r]\ar[d]&P_k\ar[d]\\
0\ar[r]&L_{l_j}\ar[r]&M_{l_j}\ar[r]\ar[d]&M_{l_{j-1}}\ar[r]\ar[d]&0\\
&&0&0
}\]
which gives rise to the exact sequence
\begin{equation}\label{seq} 0\rightarrow\Omega{M_{l_j}}\stackrel{i}{\rightarrow}\Omega{M_{l_{j-1}}}\rightarrow L_{l_j}\rightarrow 0.
\end{equation}
%Since $L_{l_j}$ is in $\Sub\Lambda_w,$ we have $\Ext^1_{\Lambda_w}(L_{l_j},\Lambda_w) = 0$ because $\Lambda_w$ is an injective object in $\Sub\Lambda_w.$
Applying $\Hom_{\Lambda_w}(\,\,,M)$ to (\ref{seq}) we get an exact sequence
$$0\rightarrow\Hom_{\Lambda_w}(L_{l_j},M)\rightarrow\Hom_{\Lambda_w}(\Omega{M_{l_{j-1}}},M)\stackrel{}{\rightarrow}\Hom_{\Lambda_w}(\Omega{M_{l_j}},M).$$ 
Fix any $g\in\Hom_{\Lambda_w}(\Omega{M_{l_j}},M).$ Since $\underline{\Hom}_{\Lambda_w}(\Omega{M_{l_j}},M)\simeq\Ext^1_{\Lambda_w}(M_{l_j},M)=0,$ 
we have that $g$ factors through a projective $\Lambda_w$-module $P$.
Since $P$ is injective in $\Sub\Lambda_w$ and $L_{l_j}\in\Sub\Lambda_w$, we have that $g$ factors through $i$ in \eqref{seq}.
\[\xymatrix@C1cm@R.3cm{
0\ar[r]&\Omega M_{l_j}\ar^i[r]\ar_g[d]\ar[dr]&\Omega M_{l_{j-1}}\ar@{.>}[d]\ar[r]&L_{l_j}\ar[r]&0\\
&M&P\ar[l]
}\]
%there is some commutative diagram 
%\[\xymatrix@R0.3cm{
%\Omega{M_{l_j}}\ar^s[rd]\ar^g[d]\\
%M&P\ar^t[l]
%}\]
%with $P$ projective. Since $\Ext^1_{\Lambda_w}(L_{l_j},P) = 0,$ there is some commuative diagram 
%\[\xymatrix@R0.3cm{
%\Omega{M_{l_j}}\ar^v[r]\ar^s[rd]&\Omega{M_{l_{j-1}}}\ar^w[d]\\
%&P
%}\]
%and hence a commutative diagram
%\[\xymatrix@R0.3cm{
%0\ar[r]&\Omega{M_{l_j}}\ar^v[r]\ar^g[d]\ar^s[dr]&\Omega{M_{l_{j-1}}}\ar[r]\ar^w[d]&L_{l_j}\ar[r]&0\\
%&M&P\ar^t[l]
%}\]
Consequently we have an exact sequence
\[0\rightarrow\Hom_{\Lambda_w}(L_{l_j},M)\rightarrow\Hom_{\Lambda_w}(\Omega{M_{l_{j-1}}},M)\rightarrow\Hom_{\Lambda_w}(\Omega{M_{l_j}},M)\rightarrow 0.\]
We know that $\Hom_{\Lambda_w}(L_{l_j},M)=\Delta_{l_j}$ by Theorem \ref{cor2_4} \textbf{(b)}.
This implies that $\Hom_{\Lambda_w}(\Omega{M_{l_{j-1}}},M)$ has a $\Delta$-filtration if $\Hom_{\Lambda_w}(\Omega{M_{l_j}},M)$ has a $\Delta$-filtration.
Using induction on the length, we have that $ U = \Hom_{\Lambda_w}(\widetilde{\Omega}M,M)$ has a $\Delta$-filtration.
Hence $ U$ is in $\mathcal{F}(\Delta).$
%and we have $\Hom_{\Lambda_w}(L_{l_1},M) = \Hom_{\Lambda_w}(M_{l_1},M) = \Delta_{l_1}.$ We have also seen that $\Hom_{\Lambda_w}(L_{l_j},M)  = \Delta_{l_j}.$ 
%We further have $\Hom_{\Lambda_w}(L_{l_t},M) = \Hom_{\Lambda_w}(\Omega{M_{l_{t-1}}},M) = \Delta_{l_t}$ (and $\Hom_{\Lambda_w}(\Omega{M_{l_t}},M)=0).$ By induction on the length of the word $w,$ we now see that all $(\Omega{M_{l_j}},M)$ have a $\Delta$-filtration. 
%We have already seen that $\Hom_{\Lambda_w}(\Delta_w,M)$ has a $\Delta$-filtration. 

\textbf{(b)(c)}\hspace{0.1cm} Since all $\Delta_j$ are in $\Hom_{\Lambda_w}(\Sub\Lambda_w,M)$ by Theorem \ref{cor2_4} \textbf{(b)}, it follows from Proposition \ref{prep2_4} that $\mathcal{F}(\Delta)\subset\Hom_{\Lambda_w}(\Sub\Lambda_w,M)\subset{}^\bot U.$
Since $ U$ is in $\mathcal{F}(\Delta),$ then $ U$ is Ext-injective in $\mathcal{F}(\Delta),$ and is hence a summand of the characteristic tilting module. Since $ U$ is already a tilting $\Gamma$-module, it must be the characteristic tilting module. 
\end{proof}

We end this section by showing that there is induced a duality between $\Sub\Lambda_w$ and $\mathcal{F}(\Delta),$ for any choice of standard cluster tilting object associated with $w.$

\begin{theorem}\label{Thm2_6}
The functors $\xymatrix{\mod\Lambda_w\ar@<0.1cm>^{F=\Hom_{\Lambda_w}(\,\,,M)}[rr]&&
\mod\Gamma\ar@<0.1cm>^{G =\Hom_\Gamma(\,\,,M)}[ll]}$
induce dualities $\xymatrix{\Sub\Lambda_w\ar@<0.1cm>^{F}[r]&
{}^\bot U = \mathcal{F}(\Delta)\ar@<0.1cm>^{G}[l]}$.
\end{theorem}

\begin{proof}
By Theorem \ref{characteristic} we have ${}^\bot U=\mathcal{F}(\Delta),$ and hence
$F(\Sub\Lambda_w) = \Hom_{\Lambda_w}(\Sub\Lambda_w,M) \subset {}^\bot U$ by Proposition \ref{prep2_4}.
Let  $Y\in\mod\Gamma.$ Then we have a surjection $\Gamma^n\rightarrow Y$ in $\mod\Gamma.$ Applying $G$ we get an injection $G(Y)\rightarrow M^n,$ showing that $G(Y)$ is in $\Sub\Lambda_w$ since $M$ is in $\Sub\Lambda_w.$ Hence we have functors
$\xymatrix{\Sub\Lambda_w\ar@<0.1cm>^F[r]&
\mathcal{F}(\Delta)\ar@<0.1cm>^G[l]}$.

We then show that $GF\simeq{\mbox{id}}$ on  $\Sub\Lambda_w.$ For $X$ in $\Sub\Lambda_w$ we have already mentioned that there is an exact  sequence $0\rightarrow X\rightarrow M_0\rightarrow M_1\rightarrow 0$ with $M_0$ and $M_1$ in $\mbox{add}M,$ and hence an exact sequence $0\rightarrow\Hom_{\Lambda_w}(M_1,M)\rightarrow\Hom_{\Lambda_w}(M_0,M)\rightarrow\Hom_{\Lambda_w}(X,M)\rightarrow 0$ in $\mod\Gamma$ since $\Ext^1_{\Lambda_w}(M_1,M)=0.$ Applying $\Hom_{\Lambda_w}(\,\,,M)$ to the last exact sequence we get an exact sequence $$0\rightarrow\Hom_{\Gamma}(\Hom_{\Lambda_w}(X,M),M)\rightarrow\Hom_{\Gamma}(\Hom_{\Lambda_w}(M_0,M),M)\rightarrow\Hom_{\Gamma}(\Hom_{\Lambda_w}(M_1,M),M).$$
When $M'$ is a summand of $M^n,$ we have an isomorphism $\Hom_\Gamma(\Hom_{\Lambda_w}(M',M),M)\simeq M',$ and we get the following commutative diagram
\[\xymatrix@R0.3cm@C0.3cm{
0\ar[r]&\Hom_\Gamma(\Hom_{\Lambda_w}(X,M),M)\ar[r]&\Hom_\Gamma(\Hom_{\Lambda_w}(M_0,M),M)\ar[r]&\Hom_\Gamma(\Hom_{\Lambda_w}(M_1,M),M)\\
0\ar[r]&X\ar[r]&M_0\ar[r]\ar^\wr[u]&M_1\ar^\wr[u]
}\]
It follows that $GF(X) = \Hom_\Gamma(\Hom_{\Lambda_w}(X,M),M)\simeq X$ as desired.

Next we show that $FG\simeq{\mbox{id}}$ on ${}^\bot U = \mathcal{F}(\Delta).$ For $\mbox{add} U$ this follows since $F(\widetilde{\Omega}M) =  U$ and hence
$$(FG) U = F(GF)(\widetilde{\Omega}M)\simeq F(\widetilde{\Omega}M)=  U.$$
Fix any $Y$ in ${}^\bot U.$
% consider the minimal left $\mbox{add} U$-approximation $Y\stackrel{f}{\rightarrow} U_0,$ which is a monomorphism, since $ U$ is a cotilting module, and hence $Y$ is in $\Sub U.$ 
%Then $\mbox{Coker}f$ is also in ${}^\bot U,$ and we take the minimal left $\mbox{add} U$-approximation $\mbox{Coker}f\rightarrow U_1,$ which is also a monomorphism. 
Since $U$ is a cotilting module, there exists an exact sequence $0\rightarrow Y\rightarrow U_0\rightarrow U_1.$ 
Applying $\Hom_\Gamma(\,\,,M)$ we get an exact sequence $\Hom_\Gamma( U_1,M)\rightarrow\Hom_\Gamma( U_0,M)\rightarrow\Hom_\Gamma(Y,M)\rightarrow 0,$ using that $M = \Hom_{\Lambda_w}(\Lambda_w,M)$ is in $\mbox{add} U=\mbox{add}\Hom_{\Lambda_w}(\widetilde{\Omega}M,M).$ 
Applying $\Hom_{\Gamma}(\,\,,M)$ we get the exact commutative diagram 
\[\xymatrix@R0.3cm@C0.3cm{
0\ar[r]&\Hom_{\Lambda_w}(\Hom_\Gamma(Y,M),M)\ar[r]&\Hom_{\Lambda_w}(\Hom_\Gamma(U_0,M),M)\ar[r]&\Hom_{\Lambda_w}(\Hom_\Gamma(U_1,M),M)\\
0\ar[r]&Y\ar[r]&U_0\ar[r]\ar^\wr[u]&U_1\ar^\wr[u]
}\]
using that $FG\simeq{\mbox{id}}$ on $\mbox{add} U.$ It follows that $(FG)Y = \Hom_{\Lambda_w}(\Hom_\Gamma(Y,M),M)\simeq Y,$ and we are done.
\end{proof}

Since, as we have seen above, $\Hom_{\Gamma}(\,\,, M):\mathcal{F}(\Delta)\rightarrow\Sub\Lambda_w$ is an exact functor, we have the following direct consequence of Theorem \ref{Thm2_6}.

\begin{corollary}
The objects in $\Sub\Lambda_w$ are exactly the objects in $\mod\Lambda_w$  which have a filtration by the layers $L_1,\ldots,L_t.$
\end{corollary}

Since the functor $G:\Sub\Lambda_w\rightarrow\mathcal{F}(\Delta)$ is not exact, it is not the case that every filtration of an object  $X$ in $\Sub\Lambda_w$ gives rise to a filtration of $F(X).$  For example, while the indecomposable projective  $\Gamma$-modules have a unique $\Delta$-composition series, there is no analogous result for the $M_i.$ As we have seen, the $M_i $ are filtered by $L_1,\ldots, L_n,$ but it can even happen that some $L_j$ is filtered by two other $L's,$ as the following example  shows.

\begin{example}
\normalfont
Let Q be the quiver 
$\xymatrix@C.3cm@R0.3cm{
{\scriptstyle 3}\ar[d]\ar[dr]\\
{\scriptstyle 2}\ar[r]&{\scriptstyle 1}}$
and $w = s_1s_2s_3s_1s_3s_2s_1,$ and $M$ the associated cluster tilting object\newline
\[\begin{smallmatrix}
1\end{smallmatrix}
\oplus
\begin{smallmatrix}
2\\
 &1\end{smallmatrix}
\oplus
\begin{smallmatrix}
 &3\\
1& &2\\
 & & &1\end{smallmatrix}
\oplus
\begin{smallmatrix}
 &1\\
2& &3\\
 &1& &2\\
 & & & &1\end{smallmatrix}
\oplus
\begin{smallmatrix}
 & &3\\
 &1& &2\\
2& & & &1\end{smallmatrix}
\oplus
\begin{smallmatrix}
 & & &2\\
 & &3& &1\\
 &1& &2& &3\\
2& & & &1& &2\\
 & & & & & & &1\end{smallmatrix}
\oplus
\begin{smallmatrix}
 & & &1\\
 & &2& &3\\
 &3& &1& &2\\
1& & & & & &1\end{smallmatrix}
\]
The $L_j$ are then the following: 
\[\begin{smallmatrix}
1\end{smallmatrix}
\ \ \ \ \ 
\begin{smallmatrix}
2\\
 &1\end{smallmatrix}
\ \ \ \ \ 
\begin{smallmatrix}
 &3\\
1& &2\\
 & & &1\end{smallmatrix}
\ \ \ \ \ 
\begin{smallmatrix}
2& &3\\
 &1& &2\\
 & & & &1\end{smallmatrix}
\ \ \ \ \ 
\begin{smallmatrix}
2\end{smallmatrix}
\ \ \ \ \ 
\begin{smallmatrix}
 & &3\\
 &1& &2& &3\\
2& & & &1& &2\\
 & & & & & & &1\end{smallmatrix}
\ \ \ \ \ 
\begin{smallmatrix}
 &3\\
1\end{smallmatrix}
\]
\vskip.5em\noindent
Then $\begin{smallmatrix}
2\\
 &1\end{smallmatrix}$
can be filtered by {\footnotesize $1,2$}, and 
$\begin{smallmatrix}
 & &2& &3\\
 &3& &1& &2\\
1& & & & & &1\end{smallmatrix}$
can be filtered by
$\begin{smallmatrix}
 &3\\
1\end{smallmatrix}$,
$\begin{smallmatrix}
2& &3\\
 &1& &2\\
 & & & &1\end{smallmatrix}$.
\end{example}

\section{Relationship between $2$-CY Frobenius categories associated with elements in Coxetr groups}\label{sc3}
 In this section we first investigate the relationship between $\Sub\Lambda_v$ and $\Sub\Lambda_w$ when $\mathbf{w} = \mathbf{u}\mathbf{v}$ is a reduced expression. 
Note that we have $\Sub\Lambda_u\subset\Sub\Lambda_w$ \cite{BIRSc}. 
We show that there is a fully faithful functor $I_u\otimes_\Lambda:\Sub\Lambda_v\rightarrow\Sub\Lambda_w,$ which preserves $\Ext^1(\,\,,\,\,),$ and which induces an  equivalence of triangulated categories between the stable category $ \underline{\Sub}\Lambda_v$ and the $2$-CY subfactor category $\mathcal{T}_{u,v}:=(\widetilde{\Omega}_{\Lambda_w}M_{\mathbf{u}})^{\bot_1}/[\widetilde{\Omega}_{\Lambda_w}M_{\mathbf{u}}]$ of $\Sub\Lambda_w.$ 
We apply this in the next section to show that $\widetilde{\Omega}M$, which gives the Ringel dual $\End_{\Lambda_w}(\widetilde{\Omega}M)$ of the quasihereditary algebra $\End_{\Lambda_w}(M)$ for a standard cluster tilting object $M$, 
belongs to the same component as $M$ in the cluster tilting graph.

Our first aim is to show that we have a fully faithful functor $I_u\otimes_\Lambda: \Sub\Lambda_v\rightarrow\Sub\Lambda_w$ which preserves $\Ext^1(\,\,,\,\,),$ as we do in the first two lemmas.

\begin{lemma}\label{L31}
\begin{enumerate}
\item We have $I_u\otimes_\Lambda\Lambda_v\simeq I_u/I_w$.
\item We have a functor $I_u\otimes_\Lambda\ :\mod\Lambda_v\rightarrow\mod\Lambda_w.$
\end{enumerate}
\end{lemma}

\begin{proof}
\textbf{(a)}\hspace{0.1cm} The exact sequence $0\rightarrow I_v\rightarrow\Lambda\rightarrow\Lambda_v\rightarrow 0$ gives rise to an exact sequence $I_u\otimes_\Lambda I_v\rightarrow I_u\otimes_\Lambda\Lambda\rightarrow I_u\otimes_\Lambda\Lambda_v\rightarrow 0.$ 
Then we have $I_u\otimes_\Lambda\Lambda_v\simeq I_u/I_u I_v = I_u/I_w$.

\textbf{(b)}\hspace{0.1cm} If $X$ is a $\Lambda_v$-module, then it is a factor module of a free $\Lambda_v$-module,
so $I_u\otimes_\Lambda X$  is a factor module of a direct sum of copies of $I_u/I_w$, which is a $\Lambda_w$-module.
Thus $I_u\otimes_\Lambda X$ is a $\Lambda_w$-module.
\end{proof}

\begin{propos}\label{L32}
We have a fully faithful functor $I_u\otimes_\Lambda: \Sub\Lambda_v\rightarrow\Sub\Lambda_w,$ which preserves $\Ext^1(\,\,,\,\,).$
\end{propos}

\begin{proof}
Without loss of generality we can assume that $Q$ is not Dynkin.

\textbf{(i)}\hspace{0.1cm} From the exact sequence of $\Lambda$-modules $0\rightarrow I_v\rightarrow\Lambda\rightarrow \Lambda_v\rightarrow 0$
we get the exact sequence $0=\mbox{Tor}_1^\Lambda(I_u,\Lambda)\rightarrow\mbox{Tor}_1^\Lambda(I_u,\Lambda_v)\rightarrow I_u\otimes_\Lambda I_v\rightarrow I_u\otimes_\Lambda\Lambda\simeq I_u.$ 
Since $I_u\otimes_\Lambda I_v = I_{u v} = I_w\subset I_u,$ using that the word $u v$ is reduced (see \cite{BIRSm}), we  see that $\mbox{Tor}_1^\Lambda(I_u,\Lambda_v) = 0.$

\textbf{(ii)}\hspace{0.1cm} For $X$ in $\Sub\Lambda_v$ we have an exact sequence
	$0\rightarrow X\rightarrow \Lambda_v^n\rightarrow Y\rightarrow 0,$ with $Y$ in $\Sub\Lambda_v$
since $\Lambda_v$ is a cotilting $\Lambda_v$-module with $\id_{\Lambda_v}\Lambda_v\le1$. Applying $I_u\otimes_\Lambda-,$ we obtain an exact sequence
$$\mbox{Tor}_2^\Lambda(I_u,Y)\rightarrow\mbox{Tor}_1^\Lambda(I_u,X)\rightarrow\mbox{Tor}_1^\Lambda(I_u,\Lambda_v^n)\rightarrow\mbox{Tor}_1^\Lambda(I_u,Y)\rightarrow I_u\otimes_\Lambda X\rightarrow I_u\otimes_\Lambda\Lambda_v^n.$$
	Since $\mbox{pd}_\Lambda I_u\leq 1,$ and we conclude that $\mbox{Tor}_1^\Lambda(I_u,X) = 0,$ and hence $\mbox{Tor}_1^\Lambda(I_u,\Sub\Lambda_v)=0.$ 
It follows that $\mbox{Tor}_1^\Lambda(I_u,Y)=0,$ and $I_u\otimes_\Lambda X$ is a submodule of $I_u\otimes_\Lambda\Lambda_v^n$.
Since $I_u\otimes_\Lambda\Lambda_v\simeq I_u/I_w$ is in $\Sub\Lambda_w,$ it follows that $I_u\otimes_\Lambda X$ is in $\Sub\Lambda_w.$ 
Hence we have a functor $I_u\otimes_\Lambda-:\Sub\Lambda_v\rightarrow\Sub\Lambda_w.$ 
%The case of Dynkin quivers can be reduced to considering the associated extended Dynkin quivers.

\textbf{(iii)}\hspace{0.1cm} We have shown that $I_u\otimes_\Lambda- = I_u\otimes_\Lambda^L-$ on $\Sub\Lambda_v$.
We know that $I_u\otimes_\Lambda^L-$ is an autoequivalence of the derived category of $\Lambda$.
Since $\Sub\Lambda_v$ and $\Sub\Lambda_w$ are extension closed full subcategories of $\mod\Lambda$ \cite{BIRSc}, we have
$$\Ext^i_{\Lambda_v}(\,\,,\,\,)=\Ext^i_{\Lambda}(\,\,,\,\,)\simeq\Ext^i_{\Lambda}(I_u\otimes_\Lambda\,\,,I_u\otimes_\Lambda\,\,)=
\Ext^i_{\Lambda_w}(I_u\otimes_\Lambda\,\,,I_u\otimes_\Lambda\,\,)$$
for $i=0,1$. Thus we have the assertion. 
\end{proof}
	
	By Proposition \ref{U is tilting} we know that $\widetilde{\Omega}_{\Lambda_w}M_{\mathbf{w}}$ is also a cluster tilting object in $\Sub\Lambda_w.$ 
For a direct summand $\widetilde{\Omega}_{\Lambda_w}M_{\mathbf{u}}$ of $\widetilde{\Omega}_{\Lambda_w}M_{\mathbf{w}}$,
%\oplus(I_n\otimes_\Lambda\Omega_{\Lambda_w}{M_{\mathbf{v}}})$ 
we consider the subfactor category $\mathcal{T}_{u,v}:=(\widetilde{\Omega}_{\Lambda_w}M_{\mathbf{u}})^{\bot_1}/[\widetilde{\Omega}_{\Lambda_w}{M_{\mathbf{u}}}]$ of $\underline{\Sub}\Lambda_w,$
and we shall show that it is triangle equivalent to $\underline{\Sub}\Lambda_v$.
We start with the following.
	
	\begin{lemma}\label{L33}
	With the previous notation we have $I_u\otimes_\Lambda\Sub\Lambda_v\subset(\widetilde{\Omega}_{\Lambda_w}M_{\mathbf{u}})^{\bot_1}.$
	\end{lemma}
\begin{proof}
The indecomposable summands of $\Omega_{\Lambda_w}M_{\mathbf{u}}$ are the indecomposable summands of $I_{u_1}/I_w$ for $\mathbf{u} = \mathbf{u}_1\mathbf{u}_2$. We have $I_{u_1}\otimes_\Lambda\Lambda_{u_2v}\simeq I_{u_1}/I_w.$
By Lemma \ref{L32} the functor $I_{u_1}\otimes_\Lambda\ : \Sub\Lambda_{u_2v}\rightarrow\Sub\Lambda_w$ preserves $\Ext^1(\,\,,\,\,),$ so we have isomorphisms 
%$$\Ext^1_{\Lambda_w}(I_{u_1}\otimes_\Lambda\Sub\Lambda_v,I_{u_1}/I_w)\simeq\Ext^1_\Lambda(I_{u_1}\otimes_\Lambda I_{u_2}\otimes_\Lambda\Sub\Lambda_v, I_{u_1}\otimes_\Lambda\Lambda/I_{u_2v})\simeq \Ext^1_{\Lambda_w}(I_{u_2}\otimes_\Lambda\Sub\Lambda_v, \Lambda/I_{u_2v}).$$
%This term is $0$ since $I_{u_2}\otimes_\Lambda\Sub\Lambda_v\subset\Sub\Lambda_{u_2v}$ and $\Lambda/I_{u_2v}$ is projective in $\Sub\Lambda_{u_2v}.$ Further 
$$\Ext^1_{\Lambda_w}(I_{u_1}/I_w,I_u\otimes_\Lambda\Sub\Lambda_v)\simeq\Ext^1_{\Lambda_w}(I_{u_1}\otimes_\Lambda\Lambda_{u_2v},I_{u_1}\otimes_\Lambda I_{u_2}\otimes_\Lambda\Sub\Lambda_v)\simeq \Ext^1_{\Lambda_{u_2v}}(\Lambda_{u_2v},I_{u_2}\otimes_\Lambda\Sub\Lambda_v)=0. $$
Thus we have the assertion.
\end{proof}

\begin{propos}\label{proof of CT}
With the above notation, we have the following:
\begin{enumerate}
\item $(I_u\otimes_\Lambda M_{\mathbf{v}})\oplus\Omega_{\Lambda_w}M_{\mathbf{u}}$ is a cluster tilting object in $\underline{\Sub}\Lambda_w.$
\item $I_u\otimes_\Lambda M_{\mathbf{v}}$ is a cluster tilting object in $ \mathcal{T}_{u,v}.$
\end{enumerate}
\end{propos}

\begin{proof}
\textbf{(a)}\hspace{0.1cm} Let $X:=(I_u\otimes_\Lambda M_{\mathbf{v}})\oplus\Omega_{\Lambda_w}M_{\mathbf{u}}$.
Since $M_{\mathbf{v}}$ is a cluster tilting object in $\underline{\Sub}\Lambda_v,$ we have $\Ext^1_{\Lambda_v}(M_{\mathbf{v}},M_{\mathbf{v}}) = 0,$ and so $\Ext^1_{\Lambda_w}(I_u\otimes_\Lambda M_{\mathbf{v}}, I_u\otimes_\Lambda M_{\mathbf{v}}) = 0$ by  Lemma \ref{L32}. 
Since $I_u\otimes_\Lambda M_{\mathbf{v}} \in (\Omega_{\Lambda_w}M_{\mathbf{u}})^{\bot_1}$ by Lemma \ref{L33}, we have 
$\Ext^1_{\Lambda_w}(\Omega_{\Lambda_w}M_{\mathbf{u}}, I_u\otimes_\Lambda M_{\mathbf{v}})=0.$ Hence $\Ext^1_{\Lambda_w}( I_u\otimes_\Lambda M_{\mathbf{v}},\Omega_{\Lambda_w}M_{\mathbf{u}})=0,$ since $\Sub\Lambda_w$ is stably $2$-CY. 
Further $\Ext^1_{\Lambda_w}(\Omega_{\Lambda_w}{M_{\mathbf{u}}},\Omega_{\Lambda_w}{M_{\mathbf{u}}})\simeq\Ext^1_{\Lambda_w}(M_{\mathbf{u}},M_{\mathbf{u}})=0,$ using that $M_{\mathbf{u}}$ is a summand of $M_{\mathbf{w}}.$ So $X$ is a rigid object in $\underline{\Sub}\Lambda_w.$ 

Let $a$ be the number of $i\in Q_0$ appearing in $\mathbf{w}$.
We know from \cite{BIRSc} that a rigid object in $\underline{\Sub}\Lambda_w$ is cluster tilting if and only if it has at least $l(w)-a$ nonisomorphic indecomposable summands.
We only have to show that the number of nonisomorphic nonprojective indecomposable summands of the $\Lambda_w$-module $X$ is at least $l(w)-a$.
Consider the following two kinds of direct summands of $X$, where $\mathbf{u}=\mathbf{u}_1\mathbf{u}_2$ and $\mathbf{v}=\mathbf{v}_1\mathbf{v}_2$ are arbitrary decompositions of words.

\textbf{(i)}\hspace{0.1cm} $\Omega_{\Lambda_w}(P_i/I_{u_1}P_i)$, where $\mathbf{u}_1$ ends at $i$ which is not the last $i$ in $\mathbf{w}$.

\textbf{(ii)}\hspace{0.1cm} $I_u\otimes_\Lambda(P_j/I_{v_1}P_j)$, where $\mathbf{v}_1$ ends at $j$ which is not the last $j$ in $\mathbf{w}$.

We will show that these $\Lambda_w$-modules are nonprojective and pairwise nonisomorphic.
Then the number of these modules is exactly $l(w)-a$, so we have that the number of nonisomorphic nonprojective indecomposable summands of the $\Lambda_w$-module $X$ is at least $l(w)-a$. This completes the proof.

Consider the module in \textbf{(i)}. Since $\Omega_{\Lambda_w}(P_i/I_{u_1}P_i)\simeq I_{u_1}P_i/I_wP_i$, this is nonprojective by the condition on $i$.
Moreover all modules in \textbf{(i)} are pairwise nonisomorphic since the functor $\Omega_{\Lambda_w}$ is an autoequivalence of $\underline{\Sub}\Lambda_w$.

Consider the module in \textbf{(ii)}. Since $I_u\otimes_\Lambda(P_j/I_{v_1}P_j)\simeq I_uP_j/I_{uv_1}P_j$, this is nonprojective by the condition on $j$.
Moreover all modules in \textbf{(ii)} are pairwise nonisomorphic since the functor $I_u\otimes_\Lambda\ $ is fully faithful by Proposition \ref{L32}.

It remains to show that the modules in \textbf{(i)} and \textbf{(ii)} are nonisomorphic. Otherwise we have
$$I_{u_1}\otimes_\Lambda(P_i/I_{u_2v}P_i)\simeq \Omega_{\Lambda_w}(P_i/I_{u_1}P_i)
\simeq I_u\otimes_\Lambda(P_j/I_{v_1}P_j)\simeq I_{u_1}\otimes_\Lambda(I_{u_2}P_j/I_{u_2v_1}P_j).$$
Since the functor $I_{u_1}\otimes_\Lambda\ $ is fully faithful by Proposition \ref{L32}, we have $P_i/I_{u_2v}P_i\simeq I_{u_2}P_j/I_{u_2v_1}P_j$.
This means that $I_{u_2}P_j/I_{u_2v_1}P_j$ is a projective $\Lambda_{u_2v}$-module.
This implies that $j$ of $\mathbf{u}_2\mathbf{v}_1$ is the last $j$ in $\mathbf{u}_2\mathbf{v}$, a contradiction to the condition in \textbf{(ii)}.

\textbf{(b)}\hspace{0.1cm} We have the assertion from \textbf{(a)} and \cite{IY}.
\end{proof}

We now prove the main result in this section.

\begin{theorem}\label{equivalence}
Let the notation be as before. Then the functor $I_u\otimes_\Lambda-: \Sub\Lambda_v\rightarrow\Sub\Lambda_w$ induces an equivalence of triangulated categories between $\underline{\Sub}\Lambda_v$ and the subfactor category $ \mathcal{T}_{u,v}$ of $\Sub\Lambda_w.$
\end{theorem}

\begin{proof}
We have seen that $M_{\mathbf{v}}$ is a cluster tilting object in $\underline{\Sub}\Lambda_v$ and $I_u\otimes_\Lambda M_{\mathbf{v}}$ is a cluster tilting object in $ \mathcal{T}_{u,v}.$ 
To show our desired equivalence, it is by \cite[4.5]{R} sufficient to show that there is induced an isomorphism $\underline{\End}_{\Lambda_v}(M_{\mathbf{v}}) \overset{\sim}{\rightarrow} \End_{ \mathcal{T}_{u,v}}(I_u\otimes_\Lambda I_v).$

By Lemma \ref{L32} there is an isomorphism $\End_{\Lambda_v}(M_{\mathbf{v}}) \overset{\sim}{\rightarrow}\End_{\Lambda_w}(I_u\otimes_\Lambda M_{\mathbf{v}})$ which induces an isomorphism
$$[\Lambda_v](M_{\mathbf{v}}) \overset{\sim}{\rightarrow}[I_u\otimes_\Lambda\Lambda_v](I_u\otimes_\Lambda M_{\mathbf{v}})
=[I_u/I_w](I_u\otimes_\Lambda M_{\mathbf{v}}),$$ 
where $[X](Y)$ is the ideal in $\End(Y)$ whose elements are the maps factoring through objects in $\mbox{add}X.$
It is sufficient to prove the equality
$$[I_u/I_w](I_u\otimes_\Lambda M_{\mathbf{v}})=[\widetilde{\Omega}_{\Lambda_w} M_{\mathbf{u}}](I_u\otimes_\Lambda M_{\mathbf{v}}).$$

Note that the indecomposable summands of $\widetilde{\Omega}_{\Lambda_w} M_{\mathbf{u}}$ are the indecomposable summands of $I_{u_1}/I_w,$ where $\mathbf{u} = \mathbf{u}_1\mathbf{u}_2.$ Thus the left side is contained in the right side.
We now show the other inclusion. The indecomposable summands of $I_u\otimes_\Lambda M_{\mathbf{v}}$ are by Lemma \ref{L31} the indecomposable summands of $I_u/I_{uv_1},$ where $\mathbf{v} = \mathbf{v}_1\mathbf{v}_2.$ 
Assume that we have a commutative diagram 
\[\xymatrix@R0.3cm@C0.3cm{
&I_{u_1}/I_w\ar^h[dr]\\
I_u/I_{u v_1}\ar^g[ur]\ar^f[rr]&&I_u/I_{u v_1}
}\]
It is sufficient to show that the image of $g$ lies in $I_u/I_w,$ or equivalently, that the composition $I_u\twoheadrightarrow I_u/I_{uv_1}\stackrel{g}{\rightarrow} I_{u_1}/I_w \subset\Lambda/I_w$ has image in $I_u/I_w.$

We apply $\Hom_\Lambda(I_u,\,\,)$ to the exact sequence $0\rightarrow I_w\rightarrow\Lambda\rightarrow\Lambda_w\rightarrow 0$ to get the exact sequence $\Hom_\Lambda(I_u,\Lambda)\rightarrow\Hom_\Lambda(I_u,\Lambda_w)\rightarrow\Ext_\Lambda^1(I_u,I_w).$ 
The last term is $0$ by \cite[III.1.13]{BIRSc}. Hence any map $I_u\rightarrow\Lambda_w$ factors through $\Lambda$.
%\[\xymatrix@R0.3cm@C0.3cm{
%&\Lambda\ar[dr]\\
%I_u\ar[ur]\ar[rr]&&\Lambda_w
%}\]
From \cite[III.1.14]{BIRSc} we know that any map $I_u\rightarrow \Lambda$ is given by the right multiplication with an element in $\Lambda.$ Thus any map  $I_u\rightarrow \Lambda_w$ has image in $I_u/I_w.$ This finishes the proof of the theorem.
\end{proof}

\section{Application to components}\label{sc4}
In this section let $w$ be an element in a Coxeter group and $\mathbf{w}$ be a reduced expression of $w$. 
Then we have cluster tilting objects $M_{\mathbf{w}}$ and $\widetilde{\Omega}_{\Lambda_w}M_{\mathbf{w}}$ in $\Sub\Lambda_w$
by Proposition \ref{U is tilting}. 
% Note that we know from Section \ref{sc2} that $\widetilde{\Omega}_{\Lambda_w}M_{\mathbf{w}}$ is also cluster tilting.
%then $\widetilde{\Omega}_{\Lambda_w}M$ lies in the same component as $M$ in the cluster tilting graph. 
Our main result here is the following.

\begin{theorem}\label{transitivity}
%Let $M_{\mathbf{w}}$ be a standard cluster tilting object in some category $\Sub\Lambda_w$ associated with a reduced word $w$ in a Coxeter group. Then 
There is a sequence of mutations of cluster tilting objects from $M_{\mathbf{w}}$ to $\widetilde{\Omega}_{\Lambda_w}M_{\mathbf{w}}$ in $\Sub\Lambda_w$ (respectively, from $M_{\mathbf{w}}$ to $\Omega_{\Lambda_w}M_{\mathbf{w}}$ in $\underline{\Sub}\Lambda_w$).
\end{theorem}

%We prove this in two steps, by establishing sequences of mutations from $(I_{i_1}\otimes_\Lambda M_{\mathbf{v}})\oplus P_{i_1}/I_wP_{i_1}$ to $M_{\mathbf{w}}$ and to $\Omega M_{\mathbf{w}}.$ 
We use induction on $l(w).$  
If $l(w) = 1,$ there is  nothing to prove. So assume $l(w)>1,$ and write $\mathbf{w}= s_{i_1}\mathbf{v},$ where $\mathbf{w}$ is a reduced expression and $s_{i_1}$ is one of the (distinguished) generators for the Coxeter group. 
Assume that the claim has been proved for reduced expressions of length less than $l(w).$
We show that there is a sequence of mutations between $\Omega_{\Lambda_w}M_{\mathbf{w}}$ and $I_{i_1}\otimes_\Lambda M_{\mathbf{v}},$ 
and between $I_{i_1}\otimes_\Lambda M_{\mathbf{v}}$ and $M_{\mathbf{w}}$ in $\underline{\Sub}\Lambda_w$.

\begin{lemma}\label{L42}
Let $\mathbf{w} = s_{i_1}\ldots s_{i_t} = s_{i_1}\mathbf{v}$ be a reduced expression. 
\begin{enumerate}
\item $I_{i_1}\otimes_\Lambda M_\mathbf{v}$ is a cluster tilting object in $\underline{\Sub}\Lambda_w$. 
\item There is a sequence of mutations of cluster tilting objects from $I_{i_1}\otimes_\Lambda M_\mathbf{v}$ to $\Omega_{\Lambda_w}M_\mathbf{w}$ in $\underline{\Sub}\Lambda_w$. 
%where $M_\mathbf{w}$ and $M_\mathbf{v}$ are the standard cluster tilting objects associated with the reduced words $\mathbf{w}$ and $\mathbf{v}.$ 
\end{enumerate}
\end{lemma}

\begin{proof}
We have $\Omega_{\Lambda_w}(P_{i_1}/I_{i_1}P_{i_1})\simeq I_{i_1}P_{i_1}/I_wP_{i_1} \simeq I_{i_1}\otimes_\Lambda P_{i_1}/I_vP_{i_1}$ and
$$\Omega_{\Lambda_w}(P_{i_r}/I_{i_1}\ldots I_{i_r}P_{i_r})\simeq I_{i_1}\ldots I_{i_r}P_{i_r}/I_wP_{i_r}\simeq I_{i_1}\otimes_\Lambda(I_{i_2}\ldots I_{i_r}P_{i_r}/I_vP_{i_r}) \simeq I_{i_1}\otimes_\Lambda\Omega_{\Lambda_v}(P_{i_r}/I_{i_2}\ldots I_{i_r}P_{i_r})$$
for $r=2,\ldots,t$. Thus we have $\Omega_{\Lambda_w}M_\mathbf{w}\simeq I_{i_1}\otimes_\Lambda\widetilde{\Omega}_{\Lambda_v}M_{\mathbf{v}}$ in $\underline{\Sub}\Lambda_w$.

By the induction assumption there is a sequence of mutations from $M_{\mathbf{v}}$ to $\Omega_{\Lambda_v}M_{\mathbf{v}}$ in $\underline{\Sub}\Lambda_v.$ 
Then by Theorem \ref{equivalence} there is a sequence of mutations from $I_{i_1}\otimes_\Lambda M_{\mathbf{v}}$ to $I_{i_1}\otimes_\Lambda\widetilde{\Omega}_{\Lambda_v}M_{\mathbf{v}}$ in $\mathcal{T}_{s_{i_1},v}$.
We have an induced sequence of mutations from $I_{i_1}\otimes_\Lambda M_{\mathbf{v}}$ to $I_{i_1}\otimes_\Lambda\widetilde{\Omega}_{\Lambda_v}M_{\mathbf{v}}$ in $\underline{\Sub}\Lambda_w$
since $\Omega_{\Lambda_w}M_{s_{i_1}}=I_{i_1}P_{i_1}/I_wP_{i_1}\simeq I_{i_1}\otimes_\Lambda(P_{i_1}/I_vP_{i_1})$ is a common direct summand of $I_{i_1}\otimes_\Lambda M_\mathbf{v}$ and $I_{i_1}\otimes_\Lambda\widetilde{\Omega}_{\Lambda_v}M_{\mathbf{v}}$ \cite{IY}.
%We have added the indecomposable projective $\Lambda_w$-module $P_{i_1}/I_{i_1}P_{i_1}$ to the first expression. 
%This is according to our discussion in Section \ref{sc3}, the only indecomposable projective $\Lambda_w$-module which may not be a summand of $I_{i_1}\otimes_\Lambda M_{\mathbf{v}}.$
%We shall now compare $\widetilde{\Omega}_{\Lambda_w}M_{\mathbf{w}}$ with $(I_{i_1}\otimes\widetilde{\Omega}_{\Lambda_v}M_{\mathbf{v}}\oplus P_{i_1}/I_wP_{i_1}).$ Let $v = s_{j_1}\ldots s_{j_t},$ and choose $r$ with $1\leq r\leq t.$ 
%So we see that the modules $\widetilde{\Omega}_{\Lambda_w}M_{\mathbf{w}}$ and $(I_{i_1}\otimes_\Lambda \widetilde{\Omega}_{\Lambda_v}M_{\mathbf{v}})\oplus P_{i_1}/I_wP_{i_1}$ coincide.
\end{proof}

In addition we have the following key step.

\begin{lemma}
There is a sequence of mutations of cluster tilting objects from $M_{\mathbf{w}}$ to $I_{i_1}\otimes_\Lambda M_{\mathbf{v}}$ in $\underline{\Sub}\Lambda_w$.
\end{lemma}

\begin{proof}
Let $1=l_1 <l_2<\ldots <l_k$ be all integers with $i:=i_{l_1} = i_{l_2}=\ldots = i_{l_k}$. 
We shall show that $\mu_{l_{k-1}}\ldots\mu_{l_1}(M_{\mathbf{w}}) \simeq I_i\otimes_\Lambda M_{\mathbf{v}}$, where $\mu_{l_j}$ denotes the mutation at the vertex $l_j.$

The summand of $I_i\otimes_\Lambda M_{\mathbf{v}}$ corresponding to some $l$ which is not one of $l_1,\ldots,l_k$ is
$$I_i\otimes_\Lambda (P_{i_l}/I_{i_2}\ldots I_{i_l}P_{i_l})\simeq I_iP_{i_l}/I_{i_1}\ldots I_{i_l}P_{i_l} = P_{i_l}/I_{i_1}\ldots I_{i_l}P_{i_l},$$
which is also a summand of $\mu_{l_{k-1}}\ldots\mu_{l_1}(M_{\mathbf{w}})$.
In the rest we shall show that the summand of $\mu_{l_{k-1}}\ldots\mu_{l_1}(M_{\mathbf{w}})$ corresponding to $l_u$ for $u=1,\ldots,k-1$ is
$$I_i\otimes_\Lambda (P_i/I_{i_2}\ldots I_{i_{l_{u+1}}}P_i)\simeq I_iP_i/I_{i_1}\ldots I_{i_{l_{u+1}}}P_i.$$

Consider the chain $P_i\supset I_{i_1}P_i\supset I_{i_1}\ldots I_{i_{l_2}}P_i\supset\ldots\supset I_{i_1}\ldots I_{i_{l_k}}P_i = I_wP_i$ of submodules of $P_i.$ 
Here we have $I_{i_1}\ldots I_{i_{l_k}}P_i = I_wP_i$ since after $i_{l_k}$ there are no vertices of type $i$.
Then we know that
$$P_i/I_{l_1}P_i\twoheadleftarrow P_i/I_{i_1}\ldots I_{i_{l_2}}P_i\twoheadleftarrow P_i/I_{i_1}\ldots I_{i_{l_3}}P_i\twoheadleftarrow\ldots\twoheadleftarrow  P_i/I_wP_i$$
is part of the quiver of $\mbox{add}M_{\mathbf{w}},$ or equivalently, the quiver of $\End_{\Lambda_w}(M_{\mathbf{w}})$ \cite{BIRSc}. 

We show that after applying $\mu_{l_{u-1}}\ldots\mu_{l_1}$ for $u = 2,\ldots,k-1$, there are exactly two arrows ending at a vertex $l_u$ which are $l_{u-1}\to l_u$ and $l_{u+1}\to l_u$.
The arrows starting or ending at $l_1,\ldots,l_k$ in the quiver of $\mbox{add}M_{\mathbf{w}}$ associated to an arrow $a$ between $i$ and $j$ are indicated in the following picture.
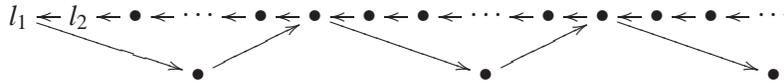
\[\xymatrix@R0.3cm@C0.3cm{
l_1\ar[drrr]&l_2\ar[l]&\bullet\ar[l]&\cdots\ar[l]&\bullet\ar[l]&\bullet\ar[l]\ar[drrr]&\bullet\ar[l]&\bullet\ar[l]&\cdots\ar[l]&\bullet\ar[l]&\bullet\ar[l]\ar[drrr]&\bullet\ar[l]&\bullet\ar[l]&\cdots\ar[l]\\
&&&\bullet\ar[urr]&&&&&\bullet\ar[urr]&&&&&\bullet
}\]
%\[\xymatrix@R0.3cm@C0.3cm{
%l_1\ar[drr]&l_2\ar[l]&\cdots\ar[l]&\bullet\ar[l]&\bullet\ar[l]\ar[drr]&\bullet\ar[l]&\cdots\ar[l]&\bullet\ar[l]&\bullet\ar[l]&\cdots\ar[l]\\
%&&\bullet\ar[urr]&&&&\bullet\ar[urr]
%}\]
(Other neighbours are omitted in this picture since the mutation behaviour is the same even if there are multiple arrows.)
%$j'$ of $i$ will be the same, so that it is sufficient to consider only one. (The argument is also the same if there are multiple arrows.)
The assertion is easily seen from performing the sequence of mutations as follows:
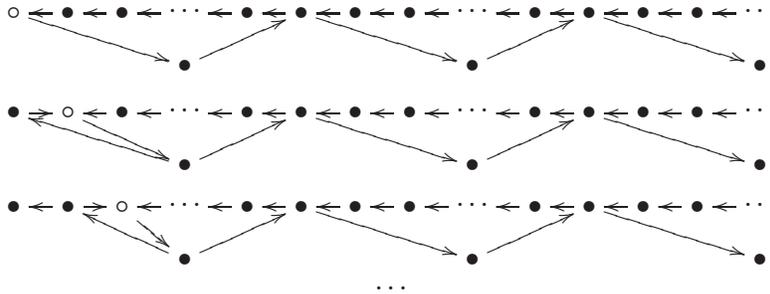
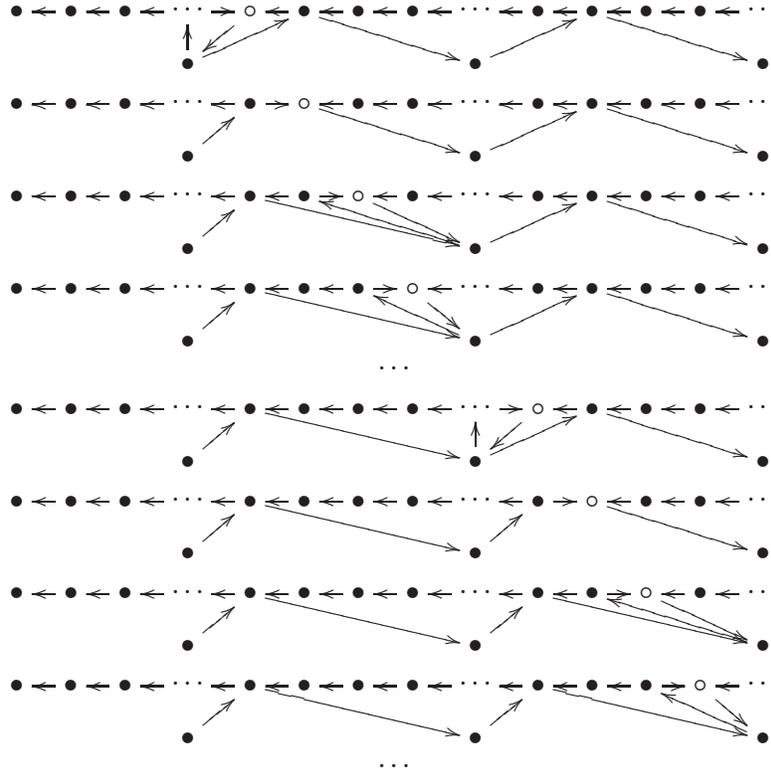
\[\xymatrix@R0.3cm@C0.3cm{
\circ\ar[drrr]&\bullet\ar[l]&\bullet\ar[l]&\cdots\ar[l]&\bullet\ar[l]&\bullet\ar[l]\ar[drrr]&\bullet\ar[l]&\bullet\ar[l]&\cdots\ar[l]&\bullet\ar[l]&\bullet\ar[l]\ar[drrr]&\bullet\ar[l]&\bullet\ar[l]&\cdots\ar[l]\\
&&&\bullet\ar[urr]&&&&&\bullet\ar[urr]&&&&&\bullet
}\]
\[\xymatrix@R0.3cm@C0.3cm{
\bullet\ar[r]&\circ\ar[drr]&\bullet\ar[l]&\cdots\ar[l]&\bullet\ar[l]&\bullet\ar[l]\ar[drrr]&\bullet\ar[l]&\bullet\ar[l]&\cdots\ar[l]&\bullet\ar[l]&\bullet\ar[l]\ar[drrr]&\bullet\ar[l]&\bullet\ar[l]&\cdots\ar[l]\\
&&&\bullet\ar[urr]\ar[ulll]&&&&&\bullet\ar[urr]&&&&&\bullet
}\]
\[\xymatrix@R0.3cm@C0.3cm{
\bullet&\bullet\ar[l]\ar[r]&\circ\ar[dr]&\cdots\ar[l]&\bullet\ar[l]&\bullet\ar[l]\ar[drrr]&\bullet\ar[l]&\bullet\ar[l]&\cdots\ar[l]&\bullet\ar[l]&\bullet\ar[l]\ar[drrr]&\bullet\ar[l]&\bullet\ar[l]&\cdots\ar[l]\\
&&&\bullet\ar[urr]\ar[ull]&&&&&\bullet\ar[urr]&&&&&\bullet
}\]
\[\cdots\]
\[\xymatrix@R0.3cm@C0.3cm{
\bullet&\bullet\ar[l]&\bullet\ar[l]&\cdots\ar[l]\ar[r]&\circ\ar[dl]&\bullet\ar[l]\ar[drrr]&\bullet\ar[l]&\bullet\ar[l]&\cdots\ar[l]&\bullet\ar[l]&\bullet\ar[l]\ar[drrr]&\bullet\ar[l]&\bullet\ar[l]&\cdots\ar[l]\\
&&&\bullet\ar[u]\ar[urr]&&&&&\bullet\ar[urr]&&&&&\bullet
}\]
\[\xymatrix@R0.3cm@C0.3cm{
\bullet&\bullet\ar[l]&\bullet\ar[l]&\cdots\ar[l]&\bullet\ar[l]\ar[r]&\circ\ar[drrr]&\bullet\ar[l]&\bullet\ar[l]&\cdots\ar[l]&\bullet\ar[l]&\bullet\ar[l]\ar[drrr]&\bullet\ar[l]&\bullet\ar[l]&\cdots\ar[l]\\
&&&\bullet\ar[ur]&&&&&\bullet\ar[urr]&&&&&\bullet
}\]
\[\xymatrix@R0.3cm@C0.3cm{
\bullet&\bullet\ar[l]&\bullet\ar[l]&\cdots\ar[l]&\bullet\ar[l]\ar[drrrr]&\bullet\ar[l]\ar[r]&\circ\ar[drr]&\bullet\ar[l]&\cdots\ar[l]&\bullet\ar[l]&\bullet\ar[l]\ar[drrr]&\bullet\ar[l]&\bullet\ar[l]&\cdots\ar[l]\\
&&&\bullet\ar[ur]&&&&&\bullet\ar[urr]\ar[ulll]&&&&&\bullet
}\]
\[\xymatrix@R0.3cm@C0.3cm{
\bullet&\bullet\ar[l]&\bullet\ar[l]&\cdots\ar[l]&\bullet\ar[l]\ar[drrrr]&\bullet\ar[l]&\bullet\ar[l]\ar[r]&\circ\ar[dr]&\cdots\ar[l]&\bullet\ar[l]&\bullet\ar[l]\ar[drrr]&\bullet\ar[l]&\bullet\ar[l]&\cdots\ar[l]\\
&&&\bullet\ar[ur]&&&&&\bullet\ar[urr]\ar[ull]&&&&&\bullet
}\]
\[\cdots\]
\[\xymatrix@R0.3cm@C0.3cm{
\bullet&\bullet\ar[l]&\bullet\ar[l]&\cdots\ar[l]&\bullet\ar[l]\ar[drrrr]&\bullet\ar[l]&\bullet\ar[l]&\bullet\ar[l]&\cdots\ar[l]\ar[r]&\circ\ar[dl]&\bullet\ar[l]\ar[drrr]&\bullet\ar[l]&\bullet\ar[l]&\cdots\ar[l]\\
&&&\bullet\ar[ur]&&&&&\bullet\ar[urr]\ar[u]&&&&&\bullet
}\]
\[\xymatrix@R0.3cm@C0.3cm{
\bullet&\bullet\ar[l]&\bullet\ar[l]&\cdots\ar[l]&\bullet\ar[l]\ar[drrrr]&\bullet\ar[l]&\bullet\ar[l]&\bullet\ar[l]&\cdots\ar[l]&\bullet\ar[l]\ar[r]&\circ\ar[drrr]&\bullet\ar[l]&\bullet\ar[l]&\cdots\ar[l]\\
&&&\bullet\ar[ur]&&&&&\bullet\ar[ur]&&&&&\bullet
}\]
\[\xymatrix@R0.3cm@C0.3cm{
\bullet&\bullet\ar[l]&\bullet\ar[l]&\cdots\ar[l]&\bullet\ar[l]\ar[drrrr]&\bullet\ar[l]&\bullet\ar[l]&\bullet\ar[l]&\cdots\ar[l]&\bullet\ar[l]\ar[drrrr]&\bullet\ar[r]\ar[l]&\circ\ar[drr]&\bullet\ar[l]&\cdots\ar[l]\\
&&&\bullet\ar[ur]&&&&&\bullet\ar[ur]&&&&&\bullet\ar[ulll]
}\]
\[\xymatrix@R0.3cm@C0.3cm{
\bullet&\bullet\ar[l]&\bullet\ar[l]&\cdots\ar[l]&\bullet\ar[l]\ar[drrrr]&\bullet\ar[l]&\bullet\ar[l]&\bullet\ar[l]&\cdots\ar[l]&\bullet\ar[l]\ar[drrrr]&\bullet\ar[l]&\bullet\ar[r]\ar[l]&\circ\ar[dr]&\cdots\ar[l]\\
&&&\bullet\ar[ur]&&&&&\bullet\ar[ur]&&&&&\bullet\ar[ull]
}\]
\[\cdots\]
%From this fact it also follows that the line with $l_j$ is kept, up to change of direction of the arrows, through the mutation process. 

%Note that arrows from vertices of type j may end at a  vertex i in the process, but only after we have made the mutation at that vertex of type i.
Under this mutation process we have replaced the indecomposable objects associated with the vertices $l_1,\ldots,l_{k-1}$ as follows. 
The object at $l_1$ is the kernel of the epimorphism $P_i/I_{i_1}\ldots I_{i_{l_2}}P_i\rightarrow P_i/I_{i_1}P_i,$ which is $I_{i_1}P_i/I_{i_1}\ldots I_{i_{l_2}}P_i.$ 
%The object at $l_2$ is the kernel of the map $I_{i_{l_1}}P_i/I_{i_{l_1}}\ldots I_{i_{l_2}}P_i\oplus P_i/I_{i_{l_1}}\ldots I_{i_{l_3}}P_i\to P_i/I_{i_{l_1}}\ldots I_{i_{l_2}}P_i,$ which is $I_{i_{l_1}}P_i/I_{i_{l_1}}\ldots I_{i_{l_3}}P_i.$ 
Similarly the object at $l_u$ is the kernel of the map
$$(I_iP_i/I_{i_1}\ldots I_{i_{l_u}}P_i)\oplus (P_i/I_{i_1}\ldots I_{i_{l_{u+1}}}P_i)\to P_i/I_{i_1}\ldots I_{i_{l_u}}P_i,$$
which is $I_iP_i/I_{i_1}\ldots I_{i_{l_{u+1}}}P_i$  for $u=1,\ldots,k-1$.
Thus we have the desired assertion.
%$I_{i_{l_1}}P_i/I_{i_{l_1}}\ldots I_{i_{l_{j+1}}}P_i$ at vertex $i_j,$ for $i<j<k,$ which is isomorphic to $I_{i_1}\otimes_\Lambda P_i/I_{i_{l_2}\ldots i_{l_{j+1}}}P_i.$
%Note that if $j\neq i,$ then $I_{i_1}\otimes_\Lambda P_j/I_{i_{l_2}\ldots i_{l_{j+1}\ldots i_t}}P_j\simeq I_{i_1}P_j/I_{i_{l_1}}\ldots I_{l_{j+1}}\ldots I_{i_t}P_j = P_j/(I_{i_{l_1}}\ldots I_{i_t})P_j,$ which is the corresponding summand of $M_{\mathbf{w}}.$
%We write \newline $P_j/(I_{i_{l_1}}\ldots I_{i_t})P_j = I_{i_1}P_j/(I_{i_{l_1}}\ldots I_{i_t})P_j \simeq I_{i_1}\otimes_\Lambda P_j/I_{i_2}\ldots I_{i_t}P_j.$ Hence also at the other vertices we get the same formula.  
\end{proof}

Putting the lemmas together we get the following main result Theorem \ref{transitivity} of this section.
%\begin{theorem}
%There is a sequence of mutations from the cluster tilting object $M_{\mathbf{w}}$ in $\Sub\Lambda_w$ to the cluster tilting object %$\widetilde{\Omega}_{\Lambda_w}M_{\mathbf{w}}.$
%\end{theorem}

We then have the following direct consequence.

\begin{corollary}
$M_{\mathbf{w}}$ and $\Omega_{\Lambda_w}M_{\mathbf{w}}$ lies in the same component of the cluster tilting graph of $\underline{\Sub}\Lambda_w$.
\end{corollary}

%\bibliographystyle{amsplain}
%\bibliography{2Aus_AlgebraRef}

\begin{thebibliography}{10}

\bibitem{AIRT}
C.~Amiot, O.~Iyama, I.~Reiten, G.~Todorov,
\emph{Preprojective algebras and c-sortable words}, arXiv:1002.4131.

\bibitem{A}
M.~Auslander, \emph{Representation dimension of artin algebras}, Queen Mary
  College Notes, 1971.

\bibitem{BKL}
M.~Barot, D.~Kussin, and H.~Lenzing, \emph{The cluster category of a canonical
  algebra}, to appear in Trans. Amer. Math. Soc., arXiv:0801.4540, 2008.

\bibitem{BIRSc}
A.~Buan, O.~Iyama, I.~Reiten, and J.~Scott, \emph{Cluster structures for
  2-{C}alabi-{Y}au categories and unipotent groups}, Compos. Math. \textbf{145} (2009), no. 4, 1035--1079.

\bibitem{BIRSm}
A.~Buan, O.~Iyama, I.~Reiten, and D.~Smith, \emph{Mutation of cluster tilting
  objects and potentials}, to appear in Amer. J. Math., arXiv:0804.3813.

\bibitem{BMRRT}
A.~Buan, R.~Marsh, M.~Reineke, I.~Reiten, and G.~Todorov, \emph{Tilting theory
  and cluster combinatorics}, Adv. Math. \textbf{204} (2006), no.~2, 572--618.

\bibitem{BMR}
A.~Buan, R.~Marsh, and I.~Reiten, \emph{Cluster-tilted algebras}, Trans. Amer.
  Math. Soc. \textbf{359} (2007), no.~1, 323--332.

\bibitem{DR}
V.~Dlab, C.~M.~Ringel, \emph{Quasi-hereditary algebras},  Illinois J. Math.  33  (1989),  no. 2, 280--291

\bibitem{GLS1}
C.~Geiss, B.~Leclere, and J.~Schr\"oer, \emph{Cluster structures and
  semicanonical bases for unipotent groups}, arXiv:math/0703039, 2007.

\bibitem{GLS2}
C.~Geiss, B.~Leclerc, and J.~Schr\"oer, \emph{Kac-{M}oody groups and cluster
  algebras}, arXiv:1001.3545, 2010.

\bibitem{I2}
O.~Iyama, \emph{Auslander correspondence}, Adv. Math. \textbf{210} (2007),
  51--82.

\bibitem{I1}
\bysame, \emph{Higher dimensional {A}uslander-{R}eiten theory on maximal
  orthogonal subcategories}, Adv. Math. \textbf{210} (2007), 22--50.

\bibitem{I3}
\bysame, \emph{Cluster tilting for higher auslander algebras},
  to appear in Adv. Math., arXiv:0809.4897, 2008.

\bibitem{IR}
O.~Iyama and I.~Reiten, \emph{{F}ormin-{Z}elevinsky muation and tilting modules
  over {C}alabi-{Y}au algebras}, Amer. J. Math. \textbf{130} (2008), no.~4,
  1087--1149.

\bibitem{IY}
O.~Iyama and Y.~Yoshino, \emph{Mutation in triangulated categories and rigid
  {C}ohen-{M}acaulay modules}, Invent. Math. \textbf{172} (2008), no.~1,
  117--168.

\bibitem{KR}
B.~Keller and I.~Reiten, \emph{Acyclic {C}alabi-{Y}au categories}, Compos.
  Math. \textbf{144} (2008), no.~5, 1332--1348, with appendix by {M}ichel {V}an
  der {B}ergh.

\bibitem{R}
C.~M.~Ringel, \emph{Iyama's finiteness theorem via strongly quasihereditary
  algebras}, JPAA (to appear).

\bibitem{R2}
\bysame, \emph{The category of modules with good filtrations over a
  quasihereditary algebra has almost split sequences}, Math. Z. \textbf{208}
  (1991), no.~2, 209--223.

\bibitem{S}
J.~Schr\"oer, \emph{Lecture at conference in Trondheim}, August 2009.

\end{thebibliography}
\providecommand{\bysame}{\leavevmode\hbox to3em{\hrulefill}\thinspace}
\providecommand{\MR}{\relax\ifhmode\unskip\space\fi MR }
\providecommand{\MRhref}[2]{%
  \href{http://www.ams.org/mathscinet-getitem?mr=#1}{#2}
}
\providecommand{\href}[2]{#2}

\end{document}